\documentstyle{amsppt}
\refstyle{A}

\nologo

\hoffset .25 true in
\voffset .2 true in

\hsize=6.1 true in
\vsize=8.5 true in

\define\imdf{\operatorname{im}d\tilde f}

\define\con{\overline{T^*_{{}_{S_\alpha}}\Cal U}}

\define\conc{\Big[ \overline{T^*_{{}_{S_\alpha}}\Cal U}\Big ]}

\define\dm{\operatorname{dim}}

\define\Adot{\bold A^\bullet}
\define\Bdot{\bold B^\bullet}
\define\Cdot{\bold C^\bullet}
\define\Fdot{\bold F^\bullet}
\define\Pdot{\bold P^\bullet}
\define\van{\phi_f\Fdot}
\define\lotimes{\ {\overset L\to \otimes}\ }

\topmatter

\title Perverse Cohomology and the Vanishing Index Theorem \endtitle

\author David B. Massey \endauthor

\address{David B. Massey, Dept. of Mathematics, Northeastern University, Boston, MA, 02115, USA} \endaddress

\email{DMASSEY\@NEU.edu}\endemail

\keywords{vanishing cycles, perverse cohomology, $a_f$ condition}\endkeywords

\subjclass{32B15, 32C35, 32C18, 32B10}\endsubjclass
\abstract
The characteristic cycle of a complex of sheaves on a complex analytic space provides weak information about the complex; essentially, it
yields the Euler characteristics of the hypercohomology of normal data to strata. We show how perverse cohomology actually allows one to
extract the individual Betti numbers of the hypercohomology of normal data to strata, not merely the Euler characteristics. We apply this
to the ``calculation'' of the vanishing cycles of a complex, and relate this to the work of Parusi\'nski and Brian\c con, Maisonobe, and
Merle on Thom's $a_f$ condition.
\endabstract
\endtopmatter

\document

\noindent\S0. {\bf Introduction}  

\vskip .2in

Let  $X$ be a $d$-dimensional complex analytic space contained in some open
subset $\Cal U$ of some
$\Bbb C^{n+1}$. Let $\tilde f:\Cal U\rightarrow\Bbb C$ be a complex analytic function and $f:=\tilde f_{|_X}$. Let $\Cal
S=\{S_\alpha\}$ be a Whitney stratification of $X$ with connected strata, and let $d_\alpha$ denote the dimension of $S_\alpha$. 

Let $R$
be a base ring which is a p.i.d., and let
$\Fdot$ be a bounded complex of sheaves of
$R$-modules which is constructible with respect to $\Cal S$; we write $\Fdot\in D^b_{{}_\Cal S}(X)$. Note
that our assumptions about
$R$ guarantee that such an $\Fdot$ is perfect (see  [{\bf G-M1}], 1.4 and [{\bf K-S}], 8.4.3). We denote the full
subcategory of
$D^b(X)$ of (complex) constructible complexes by  $D^b_{{}_\Bbb C}(X)$.

\vskip .2in

The main question which is addressed in this paper is: given the complex $\Fdot$, how does one calculate the vanishing cycle complex
$\phi_f\Fdot$? 

\vskip .2in

Of course, one should immediately ask: what does it mean to be ``given'' $\Fdot$ and what does it mean to ``calculate''
$\phi_f\Fdot$? In Theorem 2.10 of [{\bf M1}], we give an algebraic method for calculating the characteristic cycle
$\operatorname{Ch}(\phi_f\Fdot)$ when one starts with $\operatorname{Ch}(\Fdot)$ -- we refer to this result as {\it the vanishing index
theorem}. However, giving the characteristic cycle data is a far cry from describing the structure of a complex of sheaves; the
characteristic cycle supplies only Euler characteristic information about the strata. 

In Theorem 3.1 of this paper, we show how to calculate the hypercohomology of normal slices mod complex links of strata
of
$\phi_f\Fdot$, starting with analogous information about the complex $\Fdot$. While this is still fairly coarse data to associate to
$\Fdot$ and $\phi_f\Fdot$, it is certainly a substantial improvement over the characteristic cycle information.

\vskip .1in

Our technique for deriving our formulas is simple, but technical: we apply the tool of {\it perverse cohomology} (see Section 2, and [{\bf
BBD}] and [{\bf K-S}]) to the vanishing index theorem, and the Euler characteristic data almost magically turns into Betti number
data. In [{\bf M3}], we demonstrated this technique in the simple case in which the vanishing cycles are supported at an isolated
point.

\vskip .2in

Almost by definition of what we are trying to calculate, we do not need to begin with Whitney stratifications, but rather analytic
partitions for which the normal data to ``strata'' is well-defined, and then we care about only those strata which have non-trivial
hypercohomological normal data. Hence, we define {\it $\Fdot$-visible strata} and {\it $\Fdot$-normal partitionings of $X$} in 2.3 and
2.7.

\vskip .2in

In the final section of this paper, we relate the main theorem (3.1) to Thom's $a_f$ condition. In Theorem 4.4, we show that the
vanishing cycles along $f$ control Thom's $a_f$ condition. We then show how our results relate to the results of Parusi\'nski [{\bf P}]
and Brian\c con, Maisonobe, and Merle [{\bf BMM}] that Whitney stratifications adapted to $V(f)$ are $a_f$ stratifications.

\vskip .5in

\noindent\S1. {\bf The Vanishing Index Theorem}  

\vskip .2in

We continue with the notation from the introduction.

\vskip .1in

	In this section, we are going to give a general result  which describes the characteristic cycle of
$\phi_f\Fdot$ in terms of blowing-up the image of $d\tilde f$ inside the conormal spaces to strata. This is Theorem 2.10 of [{\bf
M1}].

Because $d$ is the global dimension of $X$, and we are not assuming that $X$ is pure-dimensional, or that $f$
is not constant on a
$d$-dimensional component of $X$, if $v\in\Bbb C$, then the dimension of $V(f-v)$ could be anything between $0$
and $d$. Hence, we let $\hat d_v:=1+\dm V(f-v)$, and will usually denote $\hat d_0$ by simply $\hat d$. Of
course, if we work locally, or assume that $X$ is pure-dimensional, and require
$f$ not to vanish on a component of $X$, then $\hat d$ will have attain its ``expected'' value of $d$.

\vskip .2in

For each stratum $S_\alpha\in\Cal
S$, there is a  pair $(\Bbb N_\alpha, \Bbb L_\alpha)$ consisting of a normal slice and the complex link
of $S_\alpha$; the isomorphism class of the hypercohomology $\Bbb H^*(\Bbb N_\alpha, \Bbb L_\alpha; \ \Fdot)$ is well-defined, and we
refer to it as the {\it normal data of $S_\alpha$ with respect to $\Fdot$}(see [{\bf G-M2}]).

\vskip .4in

\noindent{\bf Definition 1.1}. Recall that the {\it characteristic cycle,  $\operatorname{Ch}(\Fdot)$, of
$\Fdot$} in
$T^*\Cal U$ is the linear combination 
$\sum_{\alpha} m_{\alpha}(\Fdot) \conc$, where   the $m_\alpha(\Fdot)$ are integers determined by the Euler characteristic:  
$$m_\alpha(\Fdot)\ :=  \ (-1)^{d}\chi(\phi_{L_{|_X}}[-1]\bold F^\bullet)_{\bold x} \ = \ 
(-1)^{d}\chi(\phi_{L_{|_{{\Bbb N}_\alpha}}}[-1]{\bold F^\bullet}_{|_{{\Bbb N}_\alpha}}[-d_\alpha])_{\bold x} \ =
$$
$$(-1)^{d-d_\alpha}\chi\big(\Bbb H^*(\Bbb N_\alpha, \Bbb L_\alpha;\ \Fdot)\big)$$  for any point
$\bold x$ in 
$S_\alpha$, with normal slice ${\Bbb N}_\alpha$  at $\bold x$,  and any  $L: (\Cal U, x) \rightarrow\ (\Bbb
C,0)$ such that  $d_{\bold x}L$  is a non-degenerate covector at $\bold x$ (with respect to our fixed
stratification; see [{\bf G-M2}])  and
$L_{|_{S_\alpha}}$ has a Morse singularity at $\bold x$.  This cycle is independent of all the choices made
(see, for instance, [{\bf K-S}, Chapter IX]).

\vskip .4in

 Using the isomorphism,
$T^*\Cal U\cong
\Cal U\times\Bbb C^{n+1}$, we consider
$\operatorname{Ch}(\Fdot)$ as a cycle in
$X\times \Bbb C^{n+1}$; we use $\bold z:= (z_0, \dots, z_n)$ as coordinates on $\Cal U$ and $\bold w:= (w_0,
\dots, w_n)$ as the cotangent coordinates.  

Let $I$ denote the sheaf of ideals on $\Cal U$ given by the image of $d\tilde f$, i.e., 
$I=\big<w_0-\frac{\partial \tilde f}{\partial z_0}, \dots, w_n-\frac{\partial
\tilde f}{\partial z_n}\big>$.  For all $\alpha$, let $B_\alpha =
\operatorname{Bl}_{\operatorname{im} d\tilde f}\con$ denote the blow-up of
$\con$ along the image of $I$ in $\con$, and let $E_\alpha$ denote the corresponding exceptional divisor.  For
all
$\alpha$, we have
$E_\alpha\subseteq B_\alpha\subseteq X\times\Bbb C^{n+1}\times \Bbb P^n$.  Let $\pi: X\times\Bbb C^{n+1}\times
\Bbb P^n\rightarrow X\times \Bbb P^n$ denote the projection.  Note that, if $(\bold x, \bold w, [\eta])\in
E_\alpha$, then
$\bold w = d_\bold x\tilde f$ and so, for all
$\alpha$,
$\pi$ induces an isomorphism from $E_\alpha$ to
$\pi(E_\alpha)$.  We refer to $E := \sum_\alpha m_\alpha E_\alpha$ as the {\it total exceptional divisor} inside
the {total blow-up}
$\operatorname{Bl}_{\operatorname{im} d\tilde f}\operatorname{Ch}(\Fdot) := 
\sum_\alpha m_\alpha\operatorname{Bl}_{\operatorname{im} d\tilde f}\conc$.

\vskip .4in

Theorem 2.10 of [{\bf M1}] is:

\vskip .3in

\noindent{\bf Theorem 1.2} (Vanishing Index Theorem). {\it The projection $\pi$ induces an isomorphism between the total exceptional
divisor  
$E\subseteq\operatorname{Bl}_{\operatorname{im} d\tilde f}\operatorname{Ch}(\Fdot)$ and the sum over all
$v\in\Bbb C$ of the  (shifted) projectivized characteristic cycles of the sheaves of vanishing cycles of $\Fdot$ along
$f-v$, i.e., 
$$E\cong \pi_*(E) = \sum_{v\in\Bbb C}(-1)^{{}^{d-\hat d_v}}\Bbb P(\operatorname{Ch}(\phi_{f-v}\Fdot)).$$
 } 

\vskip .4in

\noindent{\it Remark 1.3}. Actually, there are two mild differences between the above statement and that of 2.10 of [{\bf M1}].

First of all, in [{\bf M1}], we assumed that our base ring $R$ was the complex field. However, the proof is entirely Morse-theoretic,
and hence goes through without change when $R$ is a p.i.d. (The reason for needing a p.i.d. is so that the rank of a finitely-generated
$R$-module is well-defined and is additive over long exact sequences.)

Secondly, we have introduced the factor $(-1)^{{}^{d-\hat d_v}}$; the lack of this factor was an error in [{\bf M1}]. Essentially, we
had assumed that the dimension of $V(f)$ was always $d-1$. The proof is unchanged.

\vskip .4in

We shall need one other result which follows from our work in [{\bf M1}].

\vskip .4in

\noindent{\bf Proposition 1.4}. {\it If $R$ is a p.i.d., then 
$$\big\{\bold p\in X\ |\ f(\bold p)=0, (\bold p, d_\bold p\tilde f)\in
|\operatorname{Ch}(\Fdot)|\big\}\ \subseteq\ \operatorname{supp}\phi_f\Fdot.
$$

If $R$ is a field and $\Pdot$ is a perverse sheaf on $X$, then 
$$\big\{\bold p\in X\ |\ f(\bold p)=0, (\bold p, d_\bold p\tilde f)\in
|\operatorname{Ch}(\Pdot)|\big\}\ =\ \operatorname{supp}\phi_f\Pdot.
$$
}

\vskip .2in

\noindent{\it Proof}. The first statement is immediate from Theorem 2.10 of [{\bf M1}]. The second statement follows
trivially from Theorem 3.2 of [{\bf M1}].\qed

\vskip .5in

\noindent\S2. {\bf Perverse Cohomology and Visible Strata}  

\vskip .2in

We wish to show how one can use perverse cohomology to extract Betti number information from Theorem 2.1. We  list some properties
of the perverse cohomology and of vanishing cycles that we will need later.  The reader is referred to [{\bf BBD}] and [{\bf K-S3}].

\bigskip

The perverse cohomology functor (using middle perversity, $\mu$) on $X$, 
${}^{\mu}\negmedspace H^0$, is a functor from $D^b_{{}_\Bbb C}(X)$  to the Abelian
category of perverse sheaves on $X$. One lets ${}^{\mu}\negmedspace H^i(\Fdot)$ denote ${}^{\mu}\negmedspace H^0(\Fdot[i])$.

If $\Fdot$ is constructible with respect to $\Cal S$, then ${}^{\mu}\negmedspace H^0(\Fdot)$ is also constructible with respect to
$\Cal S$, and $\big({}^{\mu}\negmedspace H^0(\Fdot)\big)_{|_{\Bbb N_\alpha}}[-d_\alpha]$ is naturally isomorphic to
${}^{\mu}\negmedspace H^0(\Fdot_{|_{\Bbb N_\alpha}}[-d_\alpha])$.

The functor 
${}^{\mu}\negmedspace H^0$, applied to a perverse sheaf $\Pdot$ is canonically isomorphic to $\Pdot$. In addition, a bounded,
constructible complex of sheaves $\Fdot$ is perverse if and only ${}^{\mu}\negmedspace H^k(\Fdot)=0$ for all $k\neq 0$. In
particular, if $X$ is a local complete intersection, then ${}^{\mu}\negmedspace H^{{\dm X}}(\Bbb Z^\bullet_X)\cong \Bbb
Z^\bullet_X[\dm X]$ and
${}^{\mu}\negmedspace H^k(\Bbb Z^\bullet_X) = 0$ if $k\neq \dm X$.

The functor 
${}^{\mu}\negmedspace H^0$ commutes with vanishing cycles with a shift of $-1$, nearby cycles with a shift of $-1$, and 
Verdier dualizing.  That is, there are natural isomorphisms
$${}^{\mu}\negmedspace H^0 \circ \phi_f[-1] \cong
\phi_f[-1] \circ {}^{\mu}\negmedspace H^0,\hskip .2in {}^{\mu}\negmedspace H^0 \circ \psi_f[-1] \cong
\psi_f[-1] \circ {}^{\mu}\negmedspace H^0, \hskip .1in\text{and }\Cal D\circ {}^{\mu}\negmedspace H^0\cong {}^{\mu}\negmedspace H^0\circ\Cal D.$$ 

Let $\bold F^\bullet$  be a bounded  complex of sheaves on $X$ which is constructible with respect to $\Cal S$. Let
$S_{\operatorname{max}}$ be a maximal stratum (i.e., one not contained in the closure of another) which is contained in the support
of $\Fdot$, and let 
$m=\dm S_{\operatorname{max}}$. Then, $\left({}^{\mu}\negmedspace H^0(\Fdot)\right)_{|_{S_{\operatorname{max}}}}$ is isomorphic (in
the derived category) to the complex which has $\left(\bold H^{-m}(\Fdot)\right)_{|_{S_{\operatorname{max}}}}$ in degree
$-m$ and zero in all other degrees. 

In particular, $\operatorname{supp}\Fdot = \bigcup_i \operatorname{supp}{}^{\mu}\negmedspace H^i(\Fdot)$, and if $\Fdot$
is supported on an isolated point, $\bold q$, then
$H^0({}^{\mu}\negmedspace H^0 (\bold F^\bullet))_{\bold q}
\cong H^0(\bold F^\bullet)_{\bold q}$. Also, $\operatorname{Ch}(\Fdot) =\sum_i(-1)^i\operatorname{Ch}\big({}^{\mu}\negmedspace
H^i(\Fdot)\big)$.

A distinguished triangle 
$$\Adot\rightarrow\Bdot\rightarrow\Cdot\rightarrow\Adot[1]$$ determines a long exact sequence in the Abelian category of perverse sheaves
$$
\dots\rightarrow {}^{\mu}\negmedspace H^{-1}(\Bdot)\rightarrow {}^{\mu}\negmedspace H^{-1}(\Cdot)\rightarrow {}^{\mu}\negmedspace
H^{0}(\Adot)\rightarrow {}^{\mu}\negmedspace H^{0}(\Bdot)\rightarrow {}^{\mu}\negmedspace H^{0}(\Cdot)\rightarrow {}^{\mu}\negmedspace
H^{1}(\Adot)\rightarrow\dots.
$$

\vskip .3in

\noindent{\bf Switching Base Rings}

\vskip .2in

In order to detect torsion by using Euler characteristics, we have to be able to switch base rings for our complexes. For each prime ideal $\frak p$ of $R$, let $k_{\frak p}$ denote the field of
fractions of $R/\frak p$, i.e., $k_0$ is the field of fractions of $R$, and for $\frak p\neq 0$, $k_{\frak p} = R/\frak p$.
There are the obvious functors $\delta_{\frak p}: \bold D^b_{{}_\Bbb C}(R_{{}_X})\rightarrow \bold D^b_{{}_\Bbb C}({(k_{\frak
p})}_{{}_X})$, which sends $\Fdot$ to $\Fdot\lotimes (k_{\frak p})^\bullet_{{}_X}$, and $\epsilon_{\frak p}: \bold
D^b_{{}_\Bbb C}({(k_{\frak p})}_{{}_X})\rightarrow\bold D^b_{{}_\Bbb C}(R_{{}_X})$, which considers $k_{\frak p}$-vector spaces as
$R$-modules. 

If $\Adot$ is a
complex of $k_{\frak p}$-vector spaces, we may consider the perverse cohomology of $\Adot$,
${}^{\mu}\negmedspace H^i_{{}_{k_{\frak p}}}(\Adot)$, or the perverse cohomology of $\epsilon(\Adot)$, which we denote by
${}^{\mu}\negmedspace H^i_{{}_R}(\Adot)$. If
$\Adot\in\bold D^b_{{}_\Bbb C}({(k_{\frak p})}_{{}_X})$ and 
$S_{\operatorname{max}}$ is a maximal stratum contained in the support of $\Adot$, then there is a canonical isomorphism
$$\epsilon\big(({}^{\mu}\negmedspace H^i_{{}_{k_{\frak p}}}(\Adot))_{|_{S_\alpha}}\big)\cong ({}^{\mu}\negmedspace
H^i_{{}_{R}}(\Adot))_{|_{S_\alpha}};$$
in particular, $\operatorname{supp}{}^{\mu}\negmedspace H^i_{{}_{k_{\frak p}}}(\Adot) = \operatorname{supp}{}^{\mu}\negmedspace
H^i_{{}_{R}}(\Adot)$.

If $\Fdot\in\bold D^b_{{}_\Bbb C}({R}_{{}_X})$, 
$S_{\operatorname{max}}$ is a maximal stratum contained in the support of $\Fdot$, and $\bold x\in S_{\operatorname{max}}$, 
then for some prime ideal
$\frak p\subset R$ and for some integer $i$, $H^i(\Fdot)_\bold x\otimes k_{\frak p}\neq 0$; it follows that
$S_{\operatorname{max}}$ is also a maximal stratum in the support of $\Fdot\lotimes {(k_{\frak p})}^\bullet_{{}_X}$.
 Thus, 
$$
\operatorname{supp}\Fdot = \bigcup_{\frak p} \operatorname{supp}(\Fdot\lotimes
{(k_{\frak p})}^\bullet_{{}_X})
$$
and so
$$
\operatorname{supp}\Fdot=\bigcup_{i, \frak p}
\operatorname{supp}{}^{\mu}\negmedspace H^i_{{}_{k_{\frak p}}}(\Fdot\lotimes {(k_{\frak p})}^\bullet_{{}_X}),
$$
where the boundedness and constructibility of $\Fdot$ imply that this union is locally finite.

For each prime ideal $\frak p$, there is a natural isomorphism in $\bold D^b_{{}_\Bbb C}(R_{{}_{V(f)}})$ given by
$$
\phi_f\big(\Fdot\lotimes{(k_{\frak p})}^\bullet_{{}_X}\big)\cong\big(\van\big)\lotimes{(k_{\frak p})}^\bullet_{{}_{V(f)}}
$$
(this is a particularly trivial case of the Sebastiani-Thom Isomorphism of [{\bf M4}]), and hence, the stalk cohomology is
given by
$$
H^i\big(\phi_f\big(\Fdot\lotimes{(k_{\frak p})}^\bullet_{{}_X}\big)\big)_\bold x\cong \big(H^i(\van)_\bold x\otimes k_{\frak
p}\big)\ \oplus\ \operatorname{Tor}\big(H^{i+1}(\van)_\bold x,\ k_{\frak p}\big).
$$

\vskip .5in

Below, we use $b_j$ to denote the $j$-th Betti number, e.g., $b_j(\Bbb N_\alpha, \Bbb
L_\alpha; \Fdot)={\operatorname{rk}}_{{}_R}\Bbb H^{j}(\Bbb N_\alpha, \Bbb
L_\alpha; \Fdot)$.

\vskip .3in

\noindent{\bf Proposition 2.1}. {\it For all integers $i$ and for all prime ideals $\frak p$ in $R$, the characteristic cycle of the
perverse cohomology of the sheaf of $\Fdot[i]$ is given by
$$
\operatorname{Ch}({}^{\mu}\negmedspace H^i(\Fdot))=(-1)^{\dm 
X}\sum_\alpha b_{i-d_\alpha}(\Bbb N_\alpha, \Bbb
L_\alpha; \Fdot)\left[\overline{T^*_{{}_{S_\alpha}}\Cal U}\right],
$$ and  the characteristic cycle of the
perverse cohomology of the sheaf of $k_{\frak p}$-vector spaces  $\Fdot[i]\lotimes (k_{\frak p})^\bullet_{{}_X}$ is given by
$$
\operatorname{Ch}({}^{\mu}\negmedspace H^i_{{}_{k_{\frak p}}}(\Fdot\lotimes (k_{\frak p})^\bullet_{{}_X}))=(-1)^{\dm 
X}\sum_\alpha c^{\frak p}_{i-d_\alpha}(\Bbb N_\alpha, \Bbb
L_\alpha; \Fdot)\left[\overline{T^*_{{}_{S_\alpha}}\Cal U}\right],
$$
where $c^{\frak p}_{j}(\Bbb N_\alpha, \Bbb
L_\alpha; \Fdot):= {\operatorname{dim}}_{k_{\frak p}}(H^{j}(\Bbb N_\alpha, \Bbb
L_\alpha; \Fdot)\otimes k_{\frak p}) + {\operatorname{dim}}_{k_{\frak
p}}\operatorname{Tor}(H^{j+1}(\Bbb N_\alpha, \Bbb
L_\alpha; \Fdot), k_{\frak p})$. }

\vskip .3in

\noindent{\it Proof}. The first formula is derived in the same way as the second, except that one has no ``$\lotimes (k_{\frak
p})^\bullet_{{}_X}$'''s anywhere; alternatively, one can deduce the first formula quickly from the second by letting $\frak p=<0>$. Thus,
we shall derive the second formula only.

\vskip .1in

Let
$\bold x\in S_\alpha$. Then,

$$m_\alpha({}^{\mu}\negmedspace H^i_{{}_{k_{\frak p}}}(\Fdot\lotimes (k_{\frak p})^\bullet_{{}_X})) =  (-1)^{\dm 
X-d_\alpha-1}\chi(\phi_{L_{|_{\Bbb N_\alpha}}}\big({{}^{\mu}\negmedspace
H^0(\Fdot[i]\lotimes (k_{\frak p})^\bullet_{{}_X})}\big)_{|_{\Bbb N_\alpha}})_{\bold x}=$$
$$(-1)^{\dm  X}\chi\big(\phi_{L_{|_{\Bbb N_\alpha}}}[-1]\big({}^{\mu}\negmedspace H^0(\Fdot[i]\lotimes (k_{\frak
p})^\bullet_{{}_X}))\big)_{|_{\Bbb N_\alpha}}[-d_\alpha]\big)_{\bold x} = $$
$$(-1)^{\dm X}\chi\big(\phi_{L_{|_{\Bbb
N_\alpha}}}[-1]{}^{\mu}\negmedspace H^0\big((\Fdot_{|_{\Bbb N_\alpha}}\lotimes (k_{\frak
p})^\bullet_{{}_{\Bbb N_\alpha}})[i-d_\alpha]\big)\big)_{\bold x}=$$
$$
 (-1)^{\dm  X}\dm_{k_{\frak p}} H^0\big(\phi_{L_{|_{\Bbb
N_\alpha}}}[-1]((\Fdot_{|_{\Bbb N_\alpha}}\lotimes (k_{\frak
p})^\bullet_{{}_{\Bbb N_\alpha}})[i-d_\alpha])\big)_{\bold x}=$$
$$
(-1)^{\dm  X}\dm_{k_{\frak p}} \Big[ \big(H^{i-d_\alpha}\big(\phi_{L_{|_{\Bbb
N_\alpha}}}[-1]\Fdot_{|_{\Bbb N_\alpha}}\big)_\bold x\otimes k_{\frak
p}\big)\ \oplus\ \operatorname{Tor}\big(H^{i-d_\alpha+1}(\phi_{L_{|_{\Bbb
N_\alpha}}}[-1]\Fdot_{|_{\Bbb N_\alpha}})_\bold x,\ k_{\frak p}\big)\Big].$$

\vskip .2in

The proof of the first statement is identical, with $k_{\frak p}$ replaced by $R$.\qed

\vskip .4in

By combining 2.1 with 1.4, we obtain

\vskip .4in

\noindent{\bf Proposition 2.2}.   {\it The following three conditions are equivalent:

\vskip .1in

\noindent a)\hskip .2in $\bold x\in\operatorname{supp}\phi_{f-f(\bold x)}\Fdot$;

\vskip .1in

\noindent b)\hskip .2inthere exists an integer
$i$ and the prime ideal
$\frak p$ in $R$ such that
$$
\bold x\in\operatorname{supp}\phi_{f-f(\bold x)}{}^{\mu}\negmedspace H^i_{{}_{k_{\frak p}}}(\Fdot\lotimes (k_{\frak
p})^\bullet_{{}_X}); 
$$
and
\vskip .1in

\noindent c)\hskip .2in there exists an integer
$i$ and the prime ideal
$\frak p$ in $R$ such that
$$
(\bold x, d_\bold x\tilde f)\in\big|\operatorname{Ch}({}^{\mu}\negmedspace H^i_{{}_{k_{\frak p}}}(\Fdot\lotimes (k_{\frak
p})^\bullet_{{}_X}))\big|.
$$ }

\vskip .3in

\noindent{\it Proof}.  By 1.4, b) and c) are equivalent. Now, a) and b) are
equivalent because
$$
\operatorname{supp}\phi_{f-f(\bold x)}\Fdot = \bigcup_{i, \frak p}\operatorname{supp}{}^{\mu}\negmedspace H^i_{{}_{k_{\frak
p}}}(\phi_{f-f(\bold x)}\Fdot\lotimes {(k_{\frak p})}^\bullet_{{}_{V(f-f(\bold x))}}) =$$ 
$$\bigcup_{i, \frak
p}\operatorname{supp}{}^{\mu}\negmedspace H^i_{{}_{k_{\frak p}}}\big(\phi_{f-f(\bold x)}(\Fdot\lotimes {(k_{\frak
p})}^\bullet_{{}_{X}})\big)=\bigcup_{i, \frak
p}\operatorname{supp}\phi_{f-f(\bold x)}{}^{\mu}\negmedspace H^i_{{}_{k_{\frak p}}}\big(\Fdot\lotimes {(k_{\frak
p})}^\bullet_{{}_{X}}\big).\qed
$$

\vskip .4in

\vskip .4in

Looking at 2.1 and 2.2, we see that any stratum $S_\alpha$ for which $H^*(\Bbb N_\alpha, \Bbb
L_\alpha;\Fdot)=0$ is essentially irrelevant as far as vanishing cycles are concerned. Hence, we make the following definition.

\vskip .3in

\noindent{\bf Definition 2.3}. A stratum $S_\alpha$ is {\it $\Fdot$-invisible} if $H^*(\Bbb N_\alpha, \Bbb
L_\alpha;\Fdot)=0$. Otherwise, we say that $S_\alpha$ is {\it $\Fdot$-visible}.

\vskip .3in

The basic principle is that $\Fdot$-invisible strata do not contribute any cohomology in Morse Theory arguments. As an example, we have 

\vskip .4in

\noindent{\bf Theorem 2.4}. {\it
$$
\bigcup_{v\in\Bbb C}\operatorname{supp}\phi_{f-v}\Fdot \ = \ \Big\{\bold x\in X\ |\ (x, d_\bold x\tilde
f)\in\bigcup\Sb\Fdot{\text-visible}\\ S_\alpha\endSb\overline{T^*_{{}_{S_\alpha}}\Cal U}\Big\}.
$$

}

\vskip .4in

\noindent{\it Proof}. This is immediate from the equivalence of a) and c) in 2.2, and the description of the characteristic cycle
given in 2.1.\qed

\vskip .4in

\noindent{\it Remark 2.5}. The union on the left side above is not just locally finite, but, in fact, locally consists of a single
support, i.e., near a point
$\bold p\in X$,
$\operatorname{supp}\phi_{f-v}\Fdot =\emptyset$ unless $v=f(\bold p)$.

Note that Theorem 2.4 immediately implies a stronger version of itself. One does not need to begin with a Whitney
stratification, but merely any stratification for which the normal data of strata with respect to $\Fdot$ is ``well-defined'', i.e.,
stratifications in which the normal data in normal slices to strata locally trivializes along the strata. Such stratifications only
require refinement by including
$\Fdot$-invisible strata in order to obtain a Whitney stratification. 

This is essentially what is required by Brian\c con, Maisonobe, and Merle in [{\bf BMM}], where they use stratifications which satisfy
Whitney's condition a) and the {\it property of local stratified triviality}.

\vskip .3in

Actually, for our purposes, we do not even need to have a ``stratification'' -- that is, we do not need the condition of the frontier.

\vskip .4in

\noindent{\bf Definition 2.6}. A {\it (complex analytic) partitioning} of $X$ is a locally finite decomposition of $X$ into disjoint
analytic submanifolds of $\Cal U$, which we still call {\it strata}, such that, for each stratum $W_\beta$, $\overline{W_\beta}$ and
$\overline{W_\beta}-W_\beta$ are closed complex analytic subsets of
$X$.

\vskip .4in

\noindent{\bf Proposition/Definition 2.7}. {\it Suppose that $\Fdot\in D^b_{{}_\Bbb C}(X)$. Let $\Cal W:=\{W_\beta\}$ be a complex
analytic partitioning of $X$ with connected strata. Then, the following  are equivalent:

\vskip .2in

\noindent a)\hskip .2in there exists a refinement $\Cal S:=\{S_\alpha\}$ of $\Cal W$ to a Whitney stratification with connected strata
such that $\Fdot\in D^b_{{}_{\Cal S}}(X)$, and such that, for all $S_\alpha$ such that
$\overline{S_\alpha}\not\in\{\overline{W_\beta}\ |\ W_\beta\in\Cal W\}$, $S_\alpha$ is $\Fdot$-invisible;

\vskip .4in

\noindent b)\hskip .2in if $\Cal S:=\{S_\alpha\}$ is a refinement of $\Cal W$ to a Whitney stratification with connected strata
such that $\Fdot\in D^b_{{}_{\Cal S}}(X)$, then, for all $S_\alpha$ such that
$\overline{S_\alpha}\not\in\{\overline{W_\beta}\ |\ W_\beta\in\Cal W\}$, $S_\alpha$ is $\Fdot$-invisible;

\vskip .4in

\noindent c)\hskip .2in for all integers $i$ and for all prime ideals $\frak p$ in $R$, 
$$
\big|\operatorname{Ch}({}^{\mu}\negmedspace H^i_{{}_{k_{\frak p}}}(\Fdot\lotimes (k_{\frak p})^\bullet_{{}_X}))\big|\subseteq
\bigcup_\beta \overline{T^*_{{}_{W_\beta}}\Cal U}.
$$

}

\vskip .2in

We refer to a partitioning $\Cal W$ which satisfies these equivalent conditions (and has connected strata) as an {\it $\Fdot$-normal
partitioning of
$X$}. Naturally, if $\Cal W$ is, in fact, a stratification of $X$ which is $\Fdot$-normal, then we refer to $\Cal W$ is an
{\it $\Fdot$-normal stratification}.

Let $\Cal W$ be an $\Fdot$-normal partitioning of $X$, and $W_\beta\in \Cal W$. Let $\Cal S$ be a Whitney refinement of $\Cal W$ such
that $\Fdot\in D^b_{{}_{\Cal S}}(X)$, and let $S_\alpha$ be the unique stratum of $\Cal S$ such that $\overline{W_\beta} =
\overline{S_\alpha}$; we refer to such an $S_\alpha$ as an {\it $\Fdot$-Whitney stratum associated to $W_\beta$}. We define the {\it
normal data of
$W_\beta$ with respect to $\Fdot$ to be the isomorphism class of the normal data of such an $S_\alpha$ with respect to $\Fdot$}, we write
$\Bbb H^*(\Bbb N_\beta, \Bbb L_\beta;\ \Fdot)$ for this normal data, and we say that
$W_\beta$ is {\it $\Fdot$-visible} if and only if it has non-zero normal data. 

{\it The definition of normal data of $W_\beta$ is independent of the
refinement $\Cal S$, and $W_\beta$ is $\Fdot$-visible if and only if there exists an integer $i$ and a prime ideal $\frak p$ in
$R$ such that 
$$
\overline{T^*_{{}_{W_\beta}}\Cal U}\ \subseteq\ \big|\operatorname{Ch}({}^{\mu}\negmedspace H^i_{{}_{k_{\frak p}}}(\Fdot\lotimes (k_{\frak
p})^\bullet_{{}_X}))\big|.
$$}

\vskip .2in

\noindent{\it Proof}. Let $\Cal S:=\{S_\alpha\}$ be a refinement of $\Cal W$ to a Whitney stratification with connected strata
such that $\Fdot\in D^b_{{}_{\Cal S}}(X)$. By 2.1, 
$$
\bigcup\Sb\Fdot{\text-visible}\\ S_\alpha\endSb\overline{T^*_{{}_{S_\alpha}}\Cal U} \ = \ \bigcup_{i, \frak
p}\big|\operatorname{Ch}({}^{\mu}\negmedspace H^i_{{}_{k_{\frak p}}}(\Fdot\lotimes (k_{\frak p})^\bullet_{{}_X}))\big|.\tag{$\dagger$}
$$

Now, 
$$\bigcup\Sb\Fdot{\text-visible}\\ S_\alpha\endSb\overline{T^*_{{}_{S_\alpha}}\Cal U}\subseteq
\bigcup_\beta\overline{T^*_{{}_{W_\beta}}\Cal U}
$$ if and only if
for all $S_\alpha$ such that
$\overline{S_\alpha}\not\in\{\overline{W_\beta}\ |\ W_\beta\in\Cal W\}$, $S_\alpha$ is $\Fdot$-invisible.

\vskip .1in

All of the conclusions follow immediately.
\qed

\vskip .5in

The point, of course, is that in results dealing with characteristic cycles, one  needs to consider only visible strata of  $\Fdot$-normal
partitions. Note that, as a result of the characterization in part c), if $\Cal W$ is an $\Fdot$-normal
partitioning of $X$, then $|\operatorname{Ch}(\Fdot)|\subseteq \bigcup_\beta \overline{T^*_{{}_{W_\beta}}\Cal U}$.

\vskip .5in

\noindent\S3. {\bf The Main Theorem}  

\vskip .2in

In this section, we are going to combine the Vanishing Index Theorem (1.2) and the characteristic cycle calculation of 2.1 to extract
the normal data of the vanishing cycles. We continue with the notations from the previous sections; especially recall the notations
used in Theorem 1.2 and in Proposition 2.1.

In what follows, if $C$ is an analytic cycle in $\Cal U$, and $V$ is a reduced and irreducible analytic subvariety of $\Cal U$, then we
write $C_V$ for the coefficient of $V$ in $C$. Also, if $Y$ and $Z$ are closed analytic subsets of $\Cal U$,  $\Cal W=\{W_\gamma\}$ is
an analytic partitioning of $Y$, and $\Omega=\sum_\gamma m_\gamma \big[\overline{T^*_{{}_{W_\gamma}}\Cal U}\big]$, then we let 
$$\Omega_{{}_{\subseteq Z}} := \sum_{W_\gamma\subseteq Z}m_\gamma\big[\overline{T^*_{{}_{W_\gamma}}\Cal U}\big],$$
and, on the set level, we define $|\Omega|_{{}_{\subseteq Z}} := |\Omega_{{}_{\subseteq Z}}|$.

\vskip .4in

\noindent{\bf Theorem 3.1}. {\it Let $\Cal S=\{S_\alpha\}$ be an $\Fdot$-normal partitioning of $X$, and let $\Cal W=\{W_\beta\}$ be a
complex analytic partitioning of $V(f)$. Let $d_\alpha:=\dm S_\alpha$ and $d_\beta:=\dm W_\beta$.

Then, $\Cal W$ is a $\phi_f\Fdot$-normal partitioning of $V(f)$ if and only if for all $\Fdot$-visible $S_\alpha$,
$$|\pi(E_\alpha)|_{{}_{\subseteq V(f)}}\subseteq\bigcup_\beta\Bbb P(\overline{T^*_{{}_{W_\beta}}\Cal U}),$$
and whenever $\Cal W$ is a $\phi_f\Fdot$-normal partitioning of $V(f)$, for all $\beta$, for all $i$, for all prime ideals $\frak p$ in
$R$,
$$
b_{i-d_\beta}(\Bbb N_\beta, \Bbb L_\beta;\ \phi_f[-1]\Fdot)\ =\ \sum_\alpha b_{i-d_\alpha}(\Bbb N_\alpha, \Bbb L_\alpha;\
\Fdot)\big[\pi(E_\alpha)\big]_{\Bbb P(\overline{T^*_{{}_{W_\beta}}\Cal U})},
$$
and
$$
c^{\frak p}_{i-d_\beta}(\Bbb N_\beta, \Bbb L_\beta;\ \phi_f[-1]\Fdot)\ =\ \sum_\alpha c^{\frak p}_{i-d_\alpha}(\Bbb N_\alpha, \Bbb
L_\alpha;\
\Fdot)\big[\pi(E_\alpha)\big]_{\Bbb P(\overline{T^*_{{}_{W_\beta}}\Cal U})}.
$$

In particular,
$$
\bigcup_{i, \frak p}\Big|\Bbb P\Big(\operatorname{Ch}\big({}^{\mu}\negmedspace H^i_{{}_{k_{\frak p}}}(\phi_f[-1]\Fdot\lotimes (k_{\frak
p})^\bullet_{{}_{V(f)}})\big)\Big)\Big|\ =\
\bigcup\Sb\Fdot{\text-visible}\\ S_\alpha\endSb\big|\pi(E_\alpha)\big|_{{}_{\subseteq V(f)}}.
$$
}

\vskip .2in

\noindent{\it Proof}. The two formulas in the theorem are derived in the same way, using Proposition 2.1. Hence, we shall prove the
notationally simpler statement about $b_*$, and conclude the statement about $c^{\frak p}_*$ at the same time.

\vskip .1in

Let $\Cal Z=\{Z_\gamma\}$ be a common Whitney refinement, with connected strata, of $\Cal S$ and $\Cal W$ such that $\Fdot$ is
constructible with respect to $\Cal Z$ and $\phi_f[-1]\Fdot$ is constructible with respect to the strata of $\Cal Z$ which are contained
in $V(f)$. Let $d_\gamma:=\dm Z_\gamma$.

 Then, the characteristic cycle of $\Pdot:={}^{\mu}\negmedspace H^i(\Fdot)$ is given by
$$
\operatorname{Ch}(\Pdot)=(-1)^d\sum_\gamma b_{i-d_\gamma}(\Bbb N_\gamma, \Bbb L_\gamma;\ \Fdot)\Big[ \overline{T^*_{{}_{Z_\gamma}}\Cal
U}\Big ].
$$
By the Vanishing Index Theorem (1.2), 
$$
(-1)^{d-\hat d}\Bbb P(\operatorname{Ch}(\phi_f\Pdot)) =(-1)^d\sum_\gamma b_{i-d_\gamma}(\Bbb N_\gamma, \Bbb L_\gamma;\
\Fdot)\big(\pi(E_\gamma)\big)_{{}_{\subseteq V(f)}}.
$$
On the other hand,
$$
(-1)^{d-\hat d}\Bbb P(\operatorname{Ch}(\phi_f\Pdot)) = (-1)^{d-\hat d}\Bbb P(\operatorname{Ch}(\phi_f{}^{\mu}\negmedspace H^i(\Fdot))) =
(-1)^{d-\hat d-1}\Bbb P(\operatorname{Ch}({}^{\mu}\negmedspace H^i(\phi_f[-1]\Fdot))) =
$$
$$
(-1)^{d}\sum_{Z_\gamma\subseteq V(f)} b_{i-d_\gamma}(\Bbb N_\gamma, \Bbb L_\gamma;\ \phi_f[-1]\Fdot)\Big[
\overline{T^*_{{}_{Z_\gamma}}\Cal U}\Big ].
$$
We conclude that if $Z_{\gamma_0}\subseteq V(f)$, then
$$
b_{i-d_{\gamma_0}}(\Bbb N_{\gamma_0}, \Bbb L_{\gamma_0};\ \phi_f[-1]\Fdot)\ =\ \sum_\gamma b_{i-d_\gamma}(\Bbb N_\gamma, \Bbb L_\gamma;\
\Fdot)\big[\pi(E_\gamma)\big]_{\Bbb P(\overline{T^*_{{}_{Z_{\gamma_0}}}\Cal U})}.
$$
In addition, since
$\Cal S$ is an $\Fdot$-normal stratification, we may replace the right-hand side by $$\sum_\alpha b_{i-d_\alpha}(\Bbb N_\alpha, \Bbb
L_\alpha;\
\Fdot)\big[\pi(E_\alpha)\big]_{\Bbb P(\overline{T^*_{{}_{Z_{\gamma_0}}}\Cal U})},$$
and we also conclude that the analogous formulas hold with $b_{i-d_\gamma}$ replaced by $c^{\frak p}_{i-d_\gamma}$. 

Using 2.7.c, it follows immediately that $\Cal W$ is a $\phi_f\Fdot$-normal partitioning of $V(f)$ if and only if for all $\Fdot$-visible
$S_\alpha$,
$$|\pi(E_\alpha)|_{{}_{\subseteq V(f)}}\subseteq\bigcup_\beta\Bbb P(\overline{T^*_{{}_{W_\beta}}\Cal U}),$$
and we are finished.\qed

\vskip .5in

\noindent{\it Remark 3.2}. The coefficient $\big[\pi(E_\alpha)\big]_{\Bbb P(\overline{T^*_{{}_{W_\beta}}\Cal U})}$ can be calculated by
moving to a generic point, $\bold x$, of $W_\beta$ and taking a normal slice $\Bbb N_\beta$. One is then reduced to the case where
$W_\beta$ consists of a single point; say $\bold 0$. By Lemma 2.7 of [{\bf M1}], the coefficient of $\Bbb
P(\overline{T^*_{{}_{0}}\Cal U})$ in $\big[\pi(E_\alpha)\big]$ can be calculated by considering the relative polar curve of
$f_{|_{S_\alpha}}$ with respect to a generic linear form $L$; the formula obtained is
$$
\big[\pi(E_\alpha)\big]_{\Bbb P(\overline{T^*_{{}_{\bold 0}}\Cal U})} = \big(\Gamma^1_{{}_{f_{|_{S_\alpha}}, L}}\cdot
V(f)\big)_\bold 0 - 
\big(\Gamma^1_{{}_{f_{|_{S_\alpha}}, L}}\cdot V(L)\big)_\bold 0.
$$

\vskip .5in

\noindent\S4. {\bf Applications to Thom's $a_f$ condition}  

\vskip .2in

In this section, we will use Theorem 3.1 to show that the vanishing cycles along $f$ control Thom's $a_f$ condition. We continue with the
notations from the previous sections.

Throughout this section, we make the simplifying assumption that $R$ is a field -- this is no real restriction, since we can always
reduce ourselves from the p.i.d. case to the field case by tensoring with $k_{\frak p}$. 

Moreover, in this section, we will usually
assume that are starting with a perverse sheaf $\Pdot$ on $X$; thus; the only possibly non-zero ${}^{\mu}\negmedspace H^i(\Pdot)$ occurs
when $i=0$. These two assumptions -- that $R$ is a field and that $\Pdot$ is perverse--greatly simplify many statements, for they imply:

\vskip .1in

\noindent a)\hskip .1in an analytic partitioning $\Cal W:=\{W_\beta\}$ of $X$ is $\Pdot$-normal if and only if
$|\operatorname{Ch}(\Pdot)|\subseteq\bigcup_\beta\overline{T^*_{{}_{W_\beta}}\Cal U}$, and

\vskip .1in

\noindent b)\hskip .1in a stratum $W_\beta$ in a $\Pdot$-normal partitioning of $X$ is $\Pdot$-visible if and only if
$\overline{T^*_{{}_{W_\beta}}\Cal U}\subseteq |\operatorname{Ch}(\Pdot)|$. 

\vskip .1in

\noindent Moreover, if $\Pdot$ is perverse, then so is
$\phi_f[-1]\Pdot$, and so statements about the vanishing cycles also simplify greatly; for example, the final statement of Theorem 3.1
becomes
$$
|\Bbb P(\operatorname{Ch}(\phi_f[-1]\Pdot))|\ =\ \bigcup\Sb\Pdot{\text-visible}\\ S_\alpha\endSb\big|\pi(E_\alpha)\big|_{{}_{\subseteq
V(f)}}.
$$
Note that, if $\bold x\in V(f)$, then there is an equality of fibres $\big(\big|\pi(E_\alpha)\big|_{{}_{\subseteq
V(f)}}\big)_\bold x= \big|\pi(E_\alpha)\big|_{\bold x}$, since, locally, the stratified critical values of $f$ are isolated.

\vskip .2in

Recall the definitions of the relative conormal space, and of Thom's $a_f$ condition in its relative conormal formulation. 

\vskip .2in

\noindent{\bf Definition 4.1}. If $M$ is an analytic submanifold of $\Cal U$ and $M\subseteq X$, then the
{\it  relative conormal space (of $M$ with respect to $f$ in $\Cal U$)}, $T^*_{f_{|_M}}\Cal U$, is given by
$$ T^*_{f_{|_M}}\Cal U:=\{(\bold x, \eta)\in T^*\Cal U\ |\ \bold x\in M,\ \eta\big(\operatorname{ker}d_\bold
x(f_{|_M})\big)=0\} = 
$$
$$
\{(\bold x, \eta)\in T^*\Cal U\ |\ \bold x\in M,\ \eta\big(T_\bold x M\cap\operatorname{ker}d_\bold x\tilde f\big)=0\}.
$$

\vskip .2in

Let $M$ and $N$ be analytic submanifolds of $X$ such that
$f$ has constant rank on $N$. Then, the pair $(M, N)$ {\it satisfies Thom's $a_f$ condition at a point $\bold
x\in N$} if and only if we have the containment
$\left(\overline{ T^*_{f_{|_M}}\Cal U}\right)_\bold x \ \subseteq \ 
\Big(T^*_{f_{|_N}}\Cal U\Big)_\bold x$ of fibres over
$\bold x$.

In particular, if $f$ is, in fact, constant on $N$, then the pair $(M, N)$ satisfies Thom's $a_f$ condition at a
point $\bold x\in N$ if and only if we have the containment
$\left(\overline{ T^*_{f_{|_M}}\Cal U}\right)_\bold x \ \subseteq \ 
\Big(T^*_{{}_N}\Cal U\Big)_\bold x$ of fibres over
$\bold x$.

\vskip .5in

\noindent{\it Remark 4.2}. We have been slightly more general in the above definition than is sometimes the case; we have not required
that the rank of $f$ be constant on
$M$. Thus, if $X$ is an analytic space, we may write that $(X_{\operatorname{reg}}, N)$ satisfies the $a_f$
condition, instead of writing the much more cumbersome $(X_{\operatorname{reg}} -
\Sigma\big(f_{|_{X_{\operatorname{reg}}}}\big), N)$ satisfies the $a_f$ condition. If $f$ is not constant on any
irreducible component of $X$, it is trivial to see that these statements are equivalent: let $\overset\circ\to X := X_{\operatorname{reg}}
-
\Sigma\big(f_{|_{X_{\operatorname{reg}}}}\big)$, which is dense in  $X_{\operatorname{reg}}$ (as $f$ is not
constant on any irreducible components of $X$), and then one shows easily that
$\overline{T^*_{f_{|_{\overset\circ\to X}}}\Cal U} =
\overline{T^*_{f_{|_{X_{\operatorname{reg}}}}}\Cal U}$.

\vskip .5in

\noindent{\bf Proposition 4.3}. {\it Suppose that $f$ is not constant on any irreducible component of $X$.   Let $E$ denote the
exceptional divisor in
$\operatorname{Bl}_{\imdf}\overline{T^*_{{}_{X_{\operatorname{reg}}}}\Cal U}\ \subseteq\ \Cal U\times\Bbb
C^{n+1}\times\Bbb P^n$. Suppose that $N\subseteq X$ is an analytic submanifold of $\Cal U$ and that $\bold x\in N$ is such that
$(X_{\operatorname{reg}}, N)$ satisfies Whitney's condition a) at $\bold x$ and such that $d_{\bold x}(f_{|_N})\equiv 0$.

Then, $(X_{\operatorname{reg}}, N)$ satisfies Whitney's $a_f$ condition at $\bold x$ if and only if there is the containment of fibres
above
$\bold x$ given by 
$$\big(\pi(E)\big)_{\bold x}\subseteq \big(\Bbb P(T^*_{{}_N}\Cal U)\big)_\bold x.$$

}

\vskip .2in

\noindent{\it Proof}. This is a modification of our proof of Theorem 4.2 of [{\bf M1}] -- we put back the Whitney a) assumption that
we went to some effort to remove in [{\bf M1}].

\vskip .1in

In Theorem 4.2 of [{\bf M1}], we showed that there is always a containment 
$$\big(\pi(E)\big)_\bold x\subseteq \Big(\Bbb P\big(\overline{T^*_{f_{|_{X_{\operatorname{reg}}}}}\Cal
U}\big)\Big)_\bold x.$$ 
It follows immediately that if $(X_{\operatorname{reg}}, N)$ satisfies Whitney's $a_f$ condition at $\bold x$, then
$\big(\pi(E)\big)_{\bold x}\subseteq \big(\Bbb P(T^*_{{}_N}\Cal U)\big)_\bold x$. We must now show the converse.

\vskip .3in

Assume that $\big(\pi(E)\big)_{\bold x}\subseteq \big(\Bbb P(T^*_{{}_N}\Cal U)\big)_\bold x$. Let $\overset\circ\to X := X_{\operatorname{reg}} - \Sigma\big(f_{|_{X_{\operatorname{reg}}}}\big)$.
 Suppose  that $[\eta]\in 
\Big(\Bbb P\big(\overline{T^*_{f_{|_{\overset\circ\to X}}}\Cal U}\big)\Big)_\bold x$. We must show that $[\eta]\in \big(\Bbb
P(T^*_{{}_N}\Cal U)\big)_\bold x$.

There exists a complex analytic
path $\alpha(t) = (\bold x(t), \eta_t)\in \overline{T^*_{f_{|_{\overset\circ\to X}}}\Cal U}$ such that
$\alpha(0) = (\bold x, \eta)$ and $\alpha(t)\in T^*_{f_{|_{\overset\circ\to X}}}\Cal U$ for $t\neq 0$. As
$f$ has no critical points on $\overset\circ\to X$, each $\eta_t$ can be written uniquely as $\eta_t =
\omega_t + \lambda(\bold x(t))d_{\bold x(t)}\tilde f$, where $\omega_t\in \big(T^*_{{}_{\overset\circ\to X}}\Cal U\big)_{\bold x(t)}$ and
$\lambda(\bold x(t))$ is a scalar. By evaluating each side on $\bold x^\prime(t)$, we find that  $\lambda(\bold x(t)) =
\frac{\eta_t(\bold x^\prime(t))}{\ \frac{d}{dt}f(\bold x(t))\ }$. 

Thus, as $\lambda(\bold x(t))$ is a quotient of two
analytic functions, there are only two possibilities for what happens to
$\lambda(\bold x(t))$ as $t\rightarrow 0$.

\vskip .3in

\noindent {\bf Case 1}:  $|\lambda(\bold x(t))|\rightarrow\infty$ as $t\rightarrow 0$.

\vskip .2in

In this case, since $\eta_t\rightarrow\eta$, it follows that 
$\dsize\frac{\eta_t}{\lambda(\bold x(t))}\rightarrow 0$ and, hence, $\dsize-\frac{\omega_t}{\lambda(\bold
x(t))}\rightarrow d_\bold x\tilde f$.   Therefore, $$\left(\bold x(t), -\frac{\omega_t}{\lambda(\bold x(t))},
\left[-\frac{\omega_t}{\lambda(\bold x(t))}-d_{\bold x(t)}\tilde f
\right]\right) \ = \ \left(\bold x(t), -\frac{\omega_t}{\lambda(\bold x(t))}, \left[\eta_t(\bold
x(t))\right]\right)\rightarrow (\bold x,
d_\bold x\tilde f, [\eta]),$$ and so $(\bold x, [\eta])\in \pi(E)$. Thus, $[\eta]\in \big(\pi(E)\big)_{\bold x}\subseteq \big(\Bbb
P(T^*_{{}_N}\Cal U)\big)_\bold x$.
 
\vskip .4in

\noindent {\bf Case 2}:  $\lambda(\bold x(t))\rightarrow\lambda_0$ as $t\rightarrow 0$. 

\vskip .2in

 In this case, $\omega_t\rightarrow\omega\in \big(\overline{T^*_{{}_{\overset\circ\to X}}\Cal U}\big)_{\bold x}$ and
$\eta=\omega+\lambda_0 d_\bold x\tilde f$. As $(X_{\operatorname{reg}}, N)$ satisfies Whitney's condition a) at $\bold x$,
$\omega\in\big(T^*_{{}_N}\Cal U\big)_\bold x$. As $d_{\bold x}(f_{|_N})\equiv 0$, $d_\bold x\tilde f\in \big(T^*_{{}_N}\Cal U\big)_\bold
x$. Thus, $\eta\in \big(T^*_{{}_N}\Cal U\big)_\bold
x$.\qed

\vskip .5in

\noindent{\bf Theorem 4.4}. {\it Let $\Pdot$ be a perverse sheaf on $X$. Let $\Cal W$ be a Whitney a) $\Pdot$-normal partitioning of $X$
such that
$V(f)$ is a union of strata and such that $\widehat{\Cal
W}:=\{W_\beta\ |\ W_\beta\subseteq V(f)\}$ is an analytic partitioning of $V(f)$. Let
$M\subseteq V(f)$ be an analytic submanifold of $\Cal U$.

Then, $\big((\overline{W_\beta})_{\operatorname{reg}}, M\big)$ satisfies the $a_f$ condition of all $\Pdot$-visible $W_\beta$ if and only
if for all
$\bold x\in M$,
$|\operatorname{Ch}(\phi_f[-1]\Pdot)|_\bold x\subseteq \big(T^*_{{}_M}\Cal U\big)_\bold x$.  }

\vskip .2in

\noindent{\it Proof}. The theorem would follow immediately from Proposition 4.3 and the fact that $$
|\Bbb P(\operatorname{Ch}(\phi_f[-1]\Pdot))|\ =\ \bigcup\Sb\Pdot{\text-visible}\\ W_\beta\endSb\big|\pi(E_\beta)\big|_{{}_{\subseteq
V(f)}},
$$
except that one must worry about strata for which $f_{|_{W_\beta}}\equiv 0$; however, in this case, $\pi(E_\beta) = \Bbb
P(\overline{T^*_{{}_{W_\beta}}\Cal U})$ (see Lemma 2.9 of [{\bf M1}]). The desired conclusion follows immediately.\qed

\vskip .5in

Thus, we see that, for a perverse sheaf, the vanishing cycles control Thom's $a_f$ condition; moreover, we could start with an
arbitrary $\Fdot\in D^b_{{}_{\Bbb C}}(X)$ and an arbitrary p.i.d. for the base ring, and then apply Theorem 4.4 to each $\Pdot:={}^{\mu}\negmedspace H^i_{{}_{k_{\frak p}}}(\Fdot\lotimes (k_{\frak
p})^\bullet_{{}_X})$.

\vskip .5in

Given the results of Parusi\'nski [{\bf P}] and Brian\c con, Maisonobe, and Merle [{\bf BMM}] that Whitney stratifications adapted to
$V(f)$ are $a_f$ stratifications (the {\it Whitney-$a_f$ result}), one might wonder if we can recover this result via the work above. The
answer is: yes and no -- basically we have to use one of the main results of [{\bf BMM}] to reach the desired conclusion. This is not
really terribly surprising; our work above has the same flavor and uses the same tools -- characteristic cycles, nearby and vanishing
cycles, perverse sheaves -- as the approach used in [{\bf BMM}]. Nonetheless, as Brian\c con, Maisonobe, and Merle concentrate on the
nearby cycles, and we wish to promote the vanishing cycles as the ``correct'' object of study, we will now indicate how to use Theorem
4.4 to conclude the Whitney-$a_f$ result.

\vskip .4in

Let $i:X-V(f)\hookrightarrow X$ and $j:V(f)\hookrightarrow X$ denote the inclusions. The following proposition follows immediately from
Theorem 3.4.2 of [{\bf BMM}] (alternatively, if follows easily from Corollary 4.6 of [{\bf M2}]).

\vskip .3in

\noindent{\bf Proposition 4.5}. {\it If $\Pdot\in D^b_{{}_\Bbb C}(X)$ is such that all of the coefficients of $\operatorname{Ch}(\Pdot)$
have the same sign (e.g., if $\Pdot$ is perverse), then
$$
|\operatorname{Ch}(i_!i^!\Pdot)|_{{}_{\subseteq V(f)}} 
= |\operatorname{Ch}(\psi_f[-1]\Pdot)|.
$$
}

\vskip .2in

\noindent{\it Proof}. Using the notation of [{\bf BMM}], the characteristic cycle of the $\Cal D$-module $M[1/f]$ is equal to
$\operatorname{Ch}(i_!i^!\Pdot)$ (see, for instance, the proof of 4.2.1 of [{\bf BMM}]). If all of the coefficients of $\operatorname{Ch}(\Pdot)$
have the same sign, there can be no cancellations in the sums appearing in 3.4.2 of [{\bf BMM}].\qed

\vskip .5in

\noindent{\bf Proposition 4.6}. {\it If $\Pdot$ is a perverse sheaf on $X$, and $\widehat{\Cal W}$ is a $j^*\Pdot$-normal
partitioning of $V(f)$, then
$$
|\operatorname{Ch}(\phi_f[-1]\Pdot)|\ \subseteq\ |\operatorname{Ch}(i_!i^!\Pdot)|_{{}_{\subseteq V(f)}} 
\ \cup\ \bigcup_{W_\beta\in \widehat{\Cal W}} \overline{T^*_{{}_{W_\beta}}\Cal U}.
$$

Moreover, if $\Cal W:=\{W_\beta\}$ is a $\Pdot$-normal
partitioning of $X$ such that $\widehat{\Cal
W}:=\{W_\beta\ |\ W_\beta\subseteq V(f)\}$ is a $j^*\Pdot$-normal partitioning
of $V(f)$, then 
$$
|\operatorname{Ch}(\phi_f[-1]\Pdot)|\ \subseteq\ \bigcup_{W_\beta\in \widehat{\Cal W}} \overline{T^*_{{}_{W_\beta}}\Cal U}.
$$
}

\vskip .2in

\noindent{\it Proof}. There is the fundamental distinguished triangle
$$j^*\Pdot[-1]\rightarrow\psi_f[-1]\Pdot\rightarrow\phi_f[-1]\Pdot\rightarrow j^*\Pdot$$ 
which  yields an equality
$$
\operatorname{Ch}(\phi_f[-1]\Pdot) \ = \ \operatorname{Ch}(\psi_f[-1]\Pdot) - \operatorname{Ch}(j^*\Pdot[-1]).
$$
Thus,
$$
|\operatorname{Ch}(\phi_f[-1]\Pdot)| \ \subseteq \ |\operatorname{Ch}(\psi_f[-1]\Pdot)|\ \cup\ |\operatorname{Ch}(j^*\Pdot[-1])|.
$$
The first statement now follows immediately from 4.5 and the fact that $\widehat{\Cal W}$ is a $j^*\Pdot$-normal
partitioning of $V(f)$.

\vskip .2in

There is a another distinguished triangle 
$$
i_!i^!\Pdot\rightarrow\Pdot\rightarrow j_*j^*\Pdot\rightarrow i_!i^!\Pdot[1],
$$
which yields the equality 
$$
\operatorname{Ch}(\Pdot) \ = \ \operatorname{Ch}(i_!i^!\Pdot) + \operatorname{Ch}(j_*j^*\Pdot).
$$
As $\operatorname{Ch}(j_*j^*\Pdot)=\pm \operatorname{Ch}(j^*\Pdot)$, we have
$|\operatorname{Ch}(j_*j^*\Pdot)|=|\operatorname{Ch}(j^*\Pdot)|$. Thus, 
$$|\operatorname{Ch}(i_!i^!\Pdot)|\subseteq
|\operatorname{Ch}(\Pdot)|\ \cup\ |\operatorname{Ch}(j^*\Pdot)|,$$
which is contained in $\bigcup_{\beta} \overline{T^*_{{}_{W_\beta}}\Cal U}$ due to our normal partitioning assumptions. Therefore,
$$|\operatorname{Ch}(i_!i^!\Pdot)|_{{}_{\subseteq V(f)}}\subseteq \bigcup_{W_\beta\in \widehat{\Cal W}} \overline{T^*_{{}_{W_\beta}}\Cal
U},$$
and the second statement now follows from the first.
\qed

\vskip .5in

\noindent{\bf Corollary 4.7}. {\it Let $\Pdot$ be a perverse sheaf on $X$. Let $\Cal W$ be a Whitney a) $\Pdot$-normal partitioning of $X$
such that
$V(f)$ is a union of strata and such that $\widehat{\Cal
W}:=\{W_\beta\ |\ W_\beta\subseteq V(f)\}$ is a $j^*\Pdot$-normal partitioning
of $V(f)$.

Let $W_\alpha\in\Cal W$ be $\Pdot$-visible and let $W_\beta\in\widehat{\Cal W}$. Then, $\big((\overline{W_\alpha})_{\operatorname{reg}},
W_\beta\big)$ satisfies the $a_f$ condition.
}

\vskip .2in

\noindent{\it Proof}. Let $\bold x\in W_\beta$. By Theorem 4.4, what we need to show is that  $|\operatorname{Ch}(\phi_f[-1]\Pdot)|_\bold
x\subseteq \big(T^*_{{}_{W_\beta}}\Cal U\big)_\bold x$. 

By 4.6, $$ |\operatorname{Ch}(\phi_f[-1]\Pdot)|\ \subseteq\ \bigcup_{W_\gamma\in
\widehat{\Cal W}} \overline{T^*_{{}_{W_\gamma}}\Cal U}.
$$
Now, Whitney's condition a) tells us that 
$$
\big(\bigcup_{W_\gamma\in \widehat{\Cal W}} \overline{T^*_{{}_{W_\gamma}}\Cal U}\big)_\bold x = \big(T^*_{{}_{W_\beta}}\Cal U\big)_\bold
x.\qed
$$

\vskip .5in

We can now give our own proof of part of the results of Parusi\'nski in [{\bf P}]
and Brian\c con, Maisonobe, and Merle in [{\bf BMM}].

\vskip .5in

\noindent{\bf Corollary 4.8}. {\it  Let $\Cal W$ be a Whitney stratification of $X$
such that
$V(f)$ is a union of strata. Let $W_\alpha, W_\beta \in\Cal W$ be such that $W_\beta\subseteq V(f)$. Then, 
$\big((\overline{W_\alpha})_{\operatorname{reg}}, W_\beta\big)$ satisfies the $a_f$ condition.  }

\vskip .2in

\noindent{\it Proof}. Let $k:\overline{W_\alpha}\hookrightarrow X$ denote the inclusion. Apply 4.7 to $\Pdot:=
k_!\,{}^{\mu}\negmedspace H^0(\Bbb C^\bullet_{{}_{\overline{W_\alpha}}}[\dm W_\alpha])\cong {}^{\mu}\negmedspace H^0(k_!\Bbb
C^\bullet_{{}_{\overline{W_\alpha}}}[\dm W_\alpha])$. 

As $k_!\Bbb
C^\bullet_{{}_{\overline{W_\alpha}}}[\dm W_\alpha]$ is constructible with respect to $\Cal W$, ${}^{\mu}\negmedspace H^0(k_!\Bbb
C^\bullet_{{}_{\overline{W_\alpha}}}[\dm W_\alpha])$ is also constructible with respect to $\Cal W$. Therefore, $\Cal W$ is a 
$\Pdot$-normal partitioning of $X$ such that
$V(f)$ is a union of strata and such that $\widehat{\Cal
W}:=\{W_\beta\ |\ W_\beta\subseteq V(f)\}$ is a $j^*\Pdot$-normal partitioning
of $V(f)$. Moreover, by construction, $W_\alpha$ is $\Pdot$-visible. The result follows.\qed

\enddocument